\input amstex
\documentstyle{amsppt}

\magnification=\magstep1
\def\PP{\bold P}

\def\mult{\operatorname{mult}}

\topmatter
\title NODAL CURVES AND POSTULATION OF GENERIC FAT POINTS ON SURFACES\endtitle
\author Edoardo BALLICO \\ Luca CHIANTINI
\endauthor
\abstract
Let $X$ be a smooth projective surface. Here we study the postulation of a general union $Z$ of fat points of $X$, when most of the connected components of $Z$ have multiplicity 2. This problem is related to the existence of "good" families of curves on $X$, with prescribed singularities, most of them being nodes, and to the cohomology of suitable line bundles on blowing ups of $X$. More precise statements are obtained in the case $X=\PP^2$.
\endabstract
\endtopmatter
\document

\heading Introduction\endheading

Let $X$ be an algebraic surface, defined over an algebraically closed field of characteristic 0; let $m>0$ be an integer and let $P\in X$. The $(m-1)$-th infinitesimal neighbourhood of $P$ in $X$ will be denoted by $mP$, hence $mP$ has $(I_{P,X})^m$ as ideal sheaf. Often $mP$ is called a {\it fat point}; $m$ is the multiplicity of $mP$ and $h^0(mP, O_{mP}) = m(m+1)/2$ is called its degree or its length. If $s,m_1,\dots,m_s$ are positive integers a nd $P_1,\dots,P_s$ are distinct points of $X$, the 0-dimensional subscheme $W =\bigcup_{1\le i\le s} m_iP_i$ of $X$ is called a {\it multi-jet} of $X$, with multiplicity $\max\{m_i\}$, type $(s;m_1,\dots,m_s)$ and length $h^0(W, O_W)$. For a fixed type $(s;m_1,\dots,m_s)$, the set of all multi-jets of type $(s;m_1,\dots,m_s)$ on $X$, is an integral variety, of dimension $2s$; hence we may speak of the {\it general} multi-jet of type $(s;m_1,\dots,m_s)$.\par\smallskip

In this paper, we study the postulation of a 0-dimensional general subscheme $Z$ of a smooth complex projective surface, under the assumption that "many" of the connected components of $Z$ are fat points of multiplicity 2. This study is a key tool for the understanding of families of curves with prescribed singularities, many of them being nodes, on a smooth surface (see e.g. [2] and [3]). This study gives also cohomological results for suitable line bundles on certain blowing ups of $X$ ([4]). Our result 0.1 below is related with the study of such families, in the blowing up of $\PP^2$ at $r$ general points.\par
We state all our results in the introduction, the proofs will be given in section 1. They use a very powerful lemma ([1], Lemma 2.3) which is a key improvement of the so-called Horace method, used in [4] for this type of problems.

\proclaim{Theorem 0.1} Fix positive integers $t,r,d_1,\dots,d_r$ and $e$; set $d_j=2$ for $r<j\le r+e$. Set $m=\max\{d_i\}_{1\le i \le r+e}$. Assume $(t+2)(t+1)/2 \ge 1 + \sum_{1\le j \le r+e}d_j(d_j+1)/2$ and $e\ge(m-1)(t-1)/2$. Then for a general multi-jet $Z:=\bigcup_{1\le j\le r+e}d_jP_j$ of type $(r+e;d_1,\dots d_{r+e})$ in $\PP^2$, we have $h^1(\PP^2,\bold I_Z(t))=0$.
\endproclaim

Here is a generalization of theorem 0.1 to the case of an arbitrary smooth projective surface.

\proclaim{Theorem 0.2} Let $X$ be a smooth projective surface. Fix integers $t>0$, $r\ge 0$, $d_j\ge 0$  $1\le j \le r$ and $e$; set $d_j=2$ for $r<j\le r+e$ and $m:=\max\{d_i\}_{1\le i \le r+e}$; assume $e\ge (m-1)(t-1)/2$. Fix $H\in $Pic($X$), with $H$ very ample and spanned. Assume $h^1(X,H^{\otimes j}) = 0$ for all $j>0$ and:\par
(1) $2(\sum_{1\le i\le r} \max\{d_i-t+j,0\})+2m\le h^0(X,H^{\otimes j})-h^0(X,H^{\otimes j-1})$ for all $2\le j\le t$\par
(2) $h^0(X,H^{\otimes j}) \ge h^0(X,H) + \sum_{1\le j \le r+e} d_j(d_j+1)/2$.\par
Then for a general multi-jet $Z:=\bigcup_{1\le j\le r+e}d_jP_j$ of type $(r+e;d_1,\dots d_{r+e})$ in $X$, we have $h^1(X, I_Z(t))=0$.
\endproclaim

Any reader of [3] will appreciate the extensions of theorem 0.1 and theorem 0.2 to the case in which we take $r$ arbitrary 0-dimensional connected subschemes, instead of $r$ multiple points (see e.g. the definition of (generalized) singularity scheme, introduced in [3], and its very effective use made there).
We will do this now.\par
Let $Z$ be a 0-dimensional connected subscheme of the germ  $\bold A_0^2$ of the affine plane at O and let $W$ be a 0-dimensional connected subscheme of a smooth projective surface $X$; set $P:=W_{red}$. We will say that $W$ is {\it equivalent} to $Z$, or that $W$ has {\it type} $Z$, if there is a formal (or \'etale, or analytic if the base field is $\bold C$) isomorphism of the germ $\bold A_0^2$ to the germ of $X$ at $P$, sending $Z$ onto $W$. The multiplicity mult${}_P(W)$ of $W$ is the maximal integer $m$ such that $Z\subset mP$. Note that mult${}_P(W)$= mult${}_O(Z)$ if $Z$ and $W$ are equivalent.\par
With these notations, the proofs of Theorems 0.1 and 0.2 give without any modification the following result:

\proclaim{Theorem 0.3} Fix positive integers $t,r,e$ and the type $Z_1,\dots,Z_r$ of $r$ 0-dimensional subschemes of the germ $\bold A^2_0$. Set $m':=\max$ $\{\mult_0(Z_i)\}_{1\le i\le r}$ and  $m:=\max\{m',2\}$. Assume
$ (t+2)(t+1)/2 \ge 1+3e+\sum_{1\le j \le r}$length$(Z_i)$
and $e\ge (m-1)(t-1)/2$. \par
Then for a general reunion $Z\subset \bold P^2$ of $e$ double points and $r$ subschemes $W_1,\dots,W_r$, with $W_i$ equivalent to $Z_i$ for every $i$, we have $h^1(\bold P^2,\bold I_Z(t))=0$.
\endproclaim

\proclaim{Theorem 0.4}Let $X$ be a smooth projective surface. Fix positive integers $t,r,e$ and the type $Z_1,\dots,Z_r$ of $r$ 0-dimensional subschemes of the germ $\bold A^2_0$. Set $m':= \max$ $\{\mult_0(Z_i)\}_{1\le i\le r}$ and  $m:=\max\{m',2\}$. Assume $e\ge (m-1)(t-1)/2$. Fix $H\in$ Pic($X$) very ample and spanned and assume $h^1(X,H^{\otimes j}) = 0$ for all $j>0$ and\par
(1) $2(\sum_{1\le i\le r} \max\{d_i-t+j,0\})+2m\le h^0(X,H^{\otimes j})-h^0(X,H^{\otimes j-1})$ for all $2\le j\le t$\par
(2) $h^0(X,H^{\otimes j}) \ge h^0(X,H) + \sum_{1\le j \le r+e} d_j(d_j+1)/2$.\par
Then for a general reunion $Z\subset X$ of $e$ double points and $r$ subschemes $W_1,\dots,W_r$, with $W_i$ equivalent to $Z_i$ for every $i$, we have $h^1(\bold P^2,\bold I_Z(t))=0$.
\endproclaim

We want to thank the referee for very useful constructive criticism on the first version of this paper. The authors were partially supported by MURST and GNSAGA of Italy.

\heading The proofs 
\endheading

We will use several times the following easy form of the so-called Horace Lemma ([4]):

\proclaim{Lemma 1.1} Let $X$ be a smooth projective surface, $H\in$ Pic($X$) an effective divisor and $Z$ a 0-dimensional subscheme of $X$. let $W:=$ Res${}_D(Z)$ be the residual scheme of $Z$ with respect to $D$, i.e. let $W\subset Z$ be the subscheme of $X$ with the conductor $(\bold I_Z:\bold I_D)$ as ideal sheaf. Set $L:= H_{|D}\in $ Pic($X$) and assume $H^1(X,\bold I_W\otimes H(-D))=H^1(D,\bold I_{Z\cap D}\otimes L)=0$. Then $H^1(X,\bold I_Z\otimes H)=0$.
\endproclaim

\demo{Proof of Theorem 0.1} If $t\le 2$ the result is trivial, hence we may assume $t\ge 3$. We have $m<t$ because otherwise $e\ge (t-1)^2/2$ and one cannot have $(t+1)(t+2)/2\ge 1+t(t+1)/2+3e$. \par
Fix a line $D\subset\PP^2$. Take a general multi-jet $W$ of type $(r;d_1,\dots,d_r)$ with length$(D\cap W)\le t+1$ and length$(D\cap W)$ as large as possible. Set $s:=t+1-$length$(D\cap W)$ and let $J$ be the union of $W$, $e-[s/2]$ general double points of $\PP^2$ and $[s/2]$ general double points supported on $D$. Note that $t\le$length$(D\cap J)\le t+1$ and that $[s/2]\le (m-1)/2$. Let $x$ be the number of connected components of $J$, with support on $D$; we have $x\ge 2$ because $m<t$ and $t+1-$length$(D\cap W)<m$, by the maximality of length$(D\cap W)$. let $m'$ be the maximum of the multiplicities of the fat points of $J\cap (\PP^2-D)$ and $e'$ be the number of double points of $J\cap(\PP^2-D)$; we have $e'\ge e-[s/2]$. If $m'<m$, since $s\le t-1$ with strict inequality when length$(D\cap J)=t$, then we have $e'\ge e-(t-1)/2>(t-1)(m-2)/2\ge(t-2)(m'-1)/2$. If $m'=m$, then $s\le (m-1)/2$, with strict inequality if length$(D\cap J)=t$; thus $e'\ge e-(s-1)/2>(t-2)(m'-1)/2.$\par
First assume length($D\cap J)=t+1$, i.e. $s$ even. By construction we have $h^0(D,\bold I_{D\cap J}(t)) =  h^0(D,\bold I_{D\cap J}(t)) = 0$. Let $G:=$ Res${}_D(J)$ be the residual scheme of $J$ with respect to $D$. By Lemma 1.1 and semicontinuity, it is sufficient to show that $h^1(\PP^2,\bold I_G(t-1)) = 0$. $G$ contains at least $e'\ge (t-2)(m'-1)/2$ double points; it is not a general multi-jet, because some of the points of its support are forced to be contained in $D$. But length$(D\cap G)$ = length ($J\cap D)-x = t+1-x\le t-1$ and we will be able to continue, exploiting again the same line, if we know how to handle the case in which length($D\cap J)=t$, i.e. $s$ is odd, for at the next step we may meet such situation.\par
Assume length($D\cap J)=t$. We take a general $P\in D$ and set $E:=J\cup\{ P \}$. Let $\bold p$ be the length 2 subscheme of $D$ with $\bold p_{red}=\{ P\}$; the scheme $\bold p$ is the second simple residue of $P$ with respect to $D$, in the sense of [1], Definition 2.2. Note that $J$ is a general multi-jet of type $(r+e-1; d_1,\dots, d_r, 2,\dots,2)$, containing $J\cap D$. Set $G':=$ Res${}_D(J)\cup\bold p$; we have $h^0(D,\bold I_{D\cap E}(t)) = h^1(D,\bold I_{D\cap E}(t)) = 0$.\par
 We claim that  by [1], Lemma 2.3, to prove 0.1 it is sufficient to prove that $h^1(\PP^2,\bold I_{G'}(t-1))=0$; since we will use the claim also to prove 0.2, 0.3 and 0.4, we want to give some details concerning the proof and translate the notations of [1], Lemma 2.1, in our situation. Set $\alpha := h^0(\PP^2,\bold I_{G'}(t-1))-$length$(G')$; the vanishing of $h^1(\PP^2, \bold I_{G'}(t-1))$ is equivalent to the fact that $\alpha\ge 0$ and that for the union, $A$ of $\alpha$ general points of $\PP^2$, we have $h^0(\PP^2,\bold I_{G'\cup A}(t-1))=0$. In the notations of the statement of [1], Lemma 2.3, we may take $Z_0=J\cup A$, $L=O_{\PP^2}(t)$, $H=D$, $r=h^0(D,O_D(t))-$length$(D\cap J)=1$ (hence the integer $r$ appearing in [1] is not our integer $r$) and $Q_1=P$, i.e. $Q_1$ is a general point of $D$; hence we obtain the claim.\par
$G'$ is not a multi-jet, but since we want to exploit again $D$ for Lemma 1.1 and $G'\cap D$ is an effective divisor on $D$ with multiplicity 2 at $P$, this is not a problem and we may repeat the construction. To obtain $H^1(\PP^2,\bold I_{g'}(t-1))=0$, we use in an essential way that  $x\ge 2$ in the following argument: since $P\in D$ and length($\bold p$) = length($\{P\})+1$, we have length($D\cap G')=2+$length(Res${}_D(J)\cap D)=2+$length($J\cap D)-x=2+t-x\le t= h^0(D,O_D(t-1))$. Alternatively, we may be sure that $J\cap D$ is not connected (i.e. that $x\ge 2$) if we impose that at each step we add at least a double point; if however at the previous step we added a double point, then at this step we are not forced to add a double point, say $2Q$ (except if $s>2$), because the residual scheme $\{Q\}=$Res${}_D(2Q)$ of $2Q$ is one connected component of $J$ and obvoiusly not the unique one, when $t\ge 2$; this alternative proof is useful for 0.2, 0.3 and 0.4.\par
To prove $h^1(\PP^2,\bold I_{G'}(t-1))=0$, we continue with the same procedure, moving some points to $D$ and taking the residue with respect to $D$; in the residue, the contribution of $\bold p$ disappears, hence we will never have more than one 0-dimensional component which is not a multiple point and this component (if any) will be a length 2 subscheme of $D$. Then we continue using the line $D$ to apply Lemma 1.1,  each time with respect to $O_{\PP^2}(t')$, with a lower integer $t'$. In this way, we finally reduce 0.1 to a maximal rank assertion for $H^0(\PP^2, O_{\PP^2}(1))$ and a 0-dimensional subscheme $A\subset \PP^2$. To conclude, it is sufficient to prove that $A$ is either empty or a reduced point. This is true because:
$$h^0(\PP^2, O_{\PP^2}(t)) \ge h^0((\PP^2, O_{\PP^2}(0))+\sum_{1\le j\le r+e} d_j(d_j+1)/2. \quad\qed$$.
\enddemo

\demo{Proof of Theorem 0.2} Fix $D\in |H|$, with $D$ smooth and irreducible. Since $H$ is very ample, we may find such $D$ passing through a general point $P$ of $X$ and tangent to an arbitrary tangent vector to $X$ at $P$. Set $L:=H_{|D}$. Since $h^1(X,H^{\otimes j})=0$ for every $j>0$, we have $h^0(X,H^{\otimes j+1})= h^1(X,H^{\otimes j})+ h^0(D,L^{\otimes j+1})$ for every $j>0$.\par
We do not want to assume the vanishing of $H^1(X,O_X)$ and this explains why, in the statement of 0.2, we are forced to add the term  $h^0(X,H)$ in equation (2).\par
the postulation of a general multi-jet on $D$ is as good as possible, i.e. for every integer $j>0$ and any datum $(x, m_1,\dots,m_x)$, then for a general multi-jet $Z$ on $D$, with datum $(x,m_1,\dots,m_x)$, the restriction map 
$H^0(D,L^{\otimes j})\to H^0(Z,L^{\otimes j}_{|Z})$ has maximal rank (see [1], Proposition 7.2). \par
We repeat verbatim the proof of 0.1. Call G(t-j) the 0-dimensional scheme that we obtain after $t-j$ steps and set $Z(j):=D\cap($Res${}_D(G(t-j))$. By the weak form of one of the assumptions in the statement of 0.2 (i.e. equation (1), without the term $2m$ in the left hand side) we have length($Z(j))\le h^0(D,L^{\otimes j-1})$  and hence the construction is possible, even if at one step we add a second residue, supported at a point of $D$. The condition on the integer length($D\cap G')$ appearing in the proof of 0.1 is satisfied because we added the term $2m$ in the left hand side of equation (1). \qed
\enddemo
\medskip

One should compare Theorem 0.1 and Theorem 0.2 with the very general paper [1], Theorem 1.1 and Corollary 1.2. After [1], the only justification for these kind of results, is given by being very explicit.

\vskip1cm
E.Ballico: Dipartimento di Matematica, Universit\'a di Trento\hskip.5cm 38050 POVO Trento (Italy)\hskip.5cm email: ballico\@science.unitn.it\par
L.Chiantini: Dipartimento di Matematica, Universit\'a di Siena, Via del Capitano 15\hskip.5cm 53100 SIENA (Italy)\hskip.5cm email: chiantini\@unisi.it

\end
\Refs\widestnumber\key{Chen}


\ref\key 1 \by Alexander J., Hirschowitz A.\paper An asymptotic vanishing theorem for generic unions of multiple points\jour preprint alg-geom 9703037 \yr 1997
\endref

\ref\key 2 \by Greuel G.M., Lossen C., Schustin E.\paper Geometry of families of nodal curves on the blown up projective plane\jour Trans. Amer. Math. Soc.\vol (to appear)
\endref

\ref\key 3 \by Greuel G.M., Lossen C., Schustin E.\paper Plane curves of minimal degree with prescribed singularities\jour Invent. Math.\vol (to appear)\endref

\ref\key 4 \by  Hirschowitz A.\paper Une conjecture pour la cohomologie des diviseurs sur les surfaces rationelles generiques\jour J. Reine Angew. Math.\vol 397 \pages 208-213\yr 1989
\endref

\endRefs
